\documentclass[11pt]{article}

\usepackage{amsxtra,amssymb,amsthm,amsmath,latexsym}
\def\R{{\mathbb R}}

\def\C{{\mathbb C}}
\def\oH{\buildrel\circ\over H}
\def\oH1{\buildrel\circ\over H\kern-.02in{}^1}

\def\l{\ell}

\newtheorem{theorem}{Theorem}[section]
\begin{document}

\title{Inverse scattering with fixed-energy data
\footnote{Math subject classification: 34R30; PACS: 03.80.+r. 03.65.Nk }
}

\author{ A.G. Ramm\\
LMA-CNRS, 31 Chemin J.Aiguier, Marseille 13402, France \\
and Mathematics Department, Kansas State University, \\
Manhattan, KS 66506-2602, USA\\
E:mail: ramm@math.ksu.edu $\quad$ 
http://www.math.ksu.edu/\,$\widetilde{\ }$\,ramm
}

\maketitle

\begin{abstract}
The Newton-Sabatier method for solving inverse scattering problem
with fixed-energy phase shifts for a sperically symmetric potential is
discussed. It is shown that this method is fundamentally wrong:
in general it cannot be carried through, the basic ansatz of R.Newton is 
wrong: the transformation kernel does not have the form postulated in 
this ansatz, in general, the method is inconsistent, and some
of the physical conclusions, e.g., existence of the transparent 
potentials, are not proved.
A mathematically justified method for solving the three-dimensional
inverse scattering problem with fixed-energy data
is described. This method is developed by A.G.Ramm for exact data and 
for noisy discrete data, and error estimates for this method are 
obtained.
Difficulties of the numerical implementation of the 
inversion method based on the 
Dirichlet-to-Neumann map are pointed out and compared with 
the difficulty of the implementation of the Ramm's inversion method.
\end{abstract}

\section{Introduction and conclusions}
 It is shown in this paper that the Newton-Sabatier (NS) method for 
solving
inverse scattering problem with fixed-energy phase shifts
for spherically-symmetric potentials is fundamentally wrong.
 The Ramm's method for solving inverse scattering problem with 
fixed-energy data is described for the exact data and for the
noisy data, and the error estimates of this method are given.
  An inversion method using the Dirichlet-to-Neumann (DN) map
is discussed, the difficulties of its numerical implementation are pointed 
out and compared with the difficulties of the implementation of the Ramm's 
inversion method.

 The inverse scattering problem (ISP) with fixed-energy data
is of basic interest in quantum and classical physics and in many
applications. The statement of the problem is well known and 
can be found in \cite{R2}.
 In \cite{N} and \cite{CS} the NS method for 
solving
ISP for spherically symmetric potentials
is described. In the
sixties P. Sabatier published several papers concerning this
procedure, and there are quite a few papers of several authors using
this procedure and generalizing it. A vast bibliography of this topic
is given in \cite{CS} and \cite{N}, and by this reason we do not
include references to many papers treating this topic.     
 
The NS method for finding $q(r), r:=|x|$, consists of the following
\cite{CS}, \cite{N}-\cite{N67}:
$$
\{\delta_{\ell}\}_{\ell=0,1,2,.....}\Rightarrow \{c_{\ell}\}\Rightarrow
K(r,s)\Rightarrow q_N(r):= -\frac{2}{r}\ \frac{d}{dr}\quad 
\frac{K(r,r)}{r}. 
\eqno{(1.1)}
$$
Here $\delta_{\ell}$ are the fixed-energy phase shifts,
generated by some potential $q(r)$, $c_{\ell}$
are some constants which should be calculated from $\delta_{\ell}$ by 
solving an infinite linear algebraic system, the $K(r,s)$ has to
be found from the equation
$$
 K(r, s) = f(r,s)-\int^r_0 K(r,t) f(t,s) t^{-2}\, dt,
\eqno{(1.2)}
$$
where 
$$
f(r,s) := \sum^\infty_{\l=0} c_\l u_\l(r) u_\l (s), \quad  u_\l (kr) := 
\sqrt{\frac{\pi kr}{2}} J_{\ell+1/2} (kr),
\eqno{(1.3)}
$$
and $J_{\ell+1/2}$ is the Bessel function.

 {\bf Our conclusion is: NS method is fundamentally wrong and, in 
general, cannot be carried through. In the exceptional cases, when it 
can be carried through, it does not yield the original generic potential 
which generated the data $\{\delta_{\ell}\}_{\ell=0,1,2,.....}$ }.

By a generic potential $q$ we mean a $q$, which is not a restriction 
to $(0,\infty)$ of an analytic function. 

It is not proved in \cite{N} and \cite{CS} that the $q_N$ from (1.1)
generates the original data $\{\delta_{\ell}\}_{\ell=0,1,2,.....}$.

{\it The R.Newton's ansatz (1.2)-(1.3) for the transformation kernel 
$K(r,s)$
is wrong: in general, $K(r,s)$ does not solve (1.2)-(1.3).}

{\it The set of potentials $q_N$ with $\sum_{\ell=0}^\infty
|c_{\ell}|<\infty$,   is not dense in   
$L_{1,1}:= \{q: \,q=\overline {q},\, ||q||<\infty, \, 
||q||:=\int_0^\infty r|q|dr\},$
 a standard scattering class. Therefore the NS method cannot be 
used even as a parameter-fitting procedure for solving ISP}.

A detailed discussion and a justification of the above conclusions
is given in \cite{R3} and \cite{R6} (see also \cite{ARS}). For convenience
of the reader we give a brief justification here.

 First, R.Newton did not prove that eq. (1.2) is solvable for all $r>0$.
If it is not solvable for at least one $r>0$, then the NS method breaks 
down: it yields $q_N$ which is not locally integrable \cite{R3}, and
the scattering theory cannot be constructed for such potentials, in 
general.

Secondly,  R. Newton did not prove  existence of the transformation kernel
$K(r,s)$, independent of $\ell$ and did not study its properties.
This was done in \cite{R4}, and it turns out that R. Newton's ansatz 
(1.2)-(1.3) for $K$, which was not justified by R. Newton in any way, is 
wrong. In particular, if $K$ solves (1.2)-(1.3), and 
$c_{\ell}$ are growing not faster than an exponential, then the function 
$K(r,r)$ must be a restriction to $(0,\infty)$ of an analytic function, 
which is not the case for a generic potential.

In the above I assumed that the inverse scattering problem, which R. 
Newton tried 
to solve, consists of finding a $q$ belonging to a certain functional 
class,
(for example, to  $L_{1,1}$, a standard scattering class), 
given the fixed-energy phase shifts generated by this $q$, and 
I have demonstrated that the NS method does not solve this problem.

In fact, the NS method does not solve even a much less interesting 
problem: given  the data
$\{\delta_{\ell}\}_{\ell=0,1,2,.....}$, corresponding to some potential 
$q$, find a potential $q_N$ which generates the same data.
This problem is much less interesting than the inverse scattering problem, 
because $q_N$ may have
no physical meaning, for example, it may deacay very slowly, etc.
 The NS method does not solve even this problem 
because equation (1.2) may be not solvable for some $r>0$,
and also because there is no proof that $q_N$ generates
the original data $\{\delta_{\ell}\}_{\ell=0,1,2,.....}$.

{\it In \cite{CS} the existence of the "transparent" potentials is 
claimed,
that is, potentials which produce, at a fixed energy, all the phase shifts 
equal to zero.} This claim is not proved because there is no proof of the 
existence of the solution of the corresponding equation (1.2) for all 
$r>0$, and also because there is no proof that $q_N$ generates
the original data $\{\delta_{\ell}\}_{\ell=0,1,2,.....}$.

Let us finally explain why the set of potentials $q_N$, which can possibly 
be obtained by the scheme (1.1) with 
$\sum_{\ell=0}^\infty |c_{\ell}|<\infty$, is not dense in $L_{1,1}$
in the norm $||\cdot||$.
Note that in the examples, considered in [4] and [5],
$c_\l=0$ for all sufficiently large $\l$, so that the condition
$\sum_{\ell=0}^\infty |c_{\ell}|<\infty$ is satisfied.

 Assume that $q\in L_{1,1}$ and that there is a 
sequence of $q_{Nj},$ denoted $q_j$ to simplify the notation, such
that $||q-q_j||\to 0$ as $j\to \infty$. Define $L(q):=\int_0^\infty rqdr$.
Clearly $L$ is a linear functional on $L_{1,1}$
and $||L||=1$. Thus $L(q_j)\to L(q)$ as $j\to \infty$. Choose $q$ with 
$L(q)\neq 0$.
Then $L(q_j)\neq 0$ for sufficiently large $j$. One has 
$K(r,r) = -\frac{r}{2} \int^r_0 sq_N(s) ds ,$ and therefore
$$K(r,r) = -\frac{L(q_j) r}{2} \left[ 1+ o(1) \right] \to \infty \hbox{\ 
as\ } r \to \infty, $$
because in our case $L(q_j)\to L(q)\neq 0$.
Thus, one has a contradiction, because 
$|K(r,r)|\leq \sup_{r>0} \sum_{\ell=0}^\infty |c_{\ell}| 
|u_{\ell}(r)||\varphi_{\ell}(r)|<\infty$,
where $\varphi_{\ell}(r)$ is the regular solution of the Schroedinger 
equation $\varphi^{\prime\prime}_\l + \varphi_\l - \frac{\l(\l +1)}{r^2}
  \varphi_\l-q(r) \varphi_\l = 0, \quad r>0$, $\varphi\sim \frac {r^\ell}
{(2\ell+1)!!}$ as $r\to 0$.
Here we have used the relations:
$$u_{\ell} (r) \sim \sin \left(r - \frac{\ell \pi}{2} \right), \quad
  \varphi_{\ell} (r) \sim |F_{\ell}|\sin \left( r-\frac{\ell \pi}{2} + 
\delta_{\ell}
\right)
  \hbox{\ as\ } r \to \infty, $$
where $\delta_{\ell}$ are the phase shifts at $k=1$,
and $F_{\ell}$ is the Jost function at $k=1$.
It can be proved that $\sup_{\ell} |F_{\ell}|<\infty$. Thus,
if $\sum^\infty_{\ell=0} |c_{\ell}| <\infty,$ then
$K(r,r) = O(1) \hbox{\ as\ } r \to \infty.$
 
In [12] a detailed analysis of the NS method is given.

In \cite{CT} a uniqueness theorem is claimed for an
equation which is a version of (1.2). In \cite{R5}
a counterexample to this wrong theorem is constructed.
In Section 2 Ramm's inversion method is 
described. In Section 3 an inversion method
which uses the DN map is described.

\section{ Ramm's inversion method for exact data}

The results we describe in this Section are taken from \cite{R2} and 
\cite{R7}. 
Assume $q\in Q:=Q_a \cap L^\infty (\R^3),$ where
$Q_a:=\{ q: q(x) = \overline{q(x)}, \quad
        q(x) \in L^2(B_a), \quad q(x) = 0
        \hbox{\ if\ } |x|\geq a \},$ $B_a:=\{x: |x|\leq a\}$.
Let $A(\alpha^\prime, \alpha)$ be the corresponding scattering amplitude 
at a fixed energy $k^2$, $k=1$ is taken without loss of generality. One 
has:
$$
A(\alpha^\prime, \alpha) = \sum^\infty_{\l =0} A_\l (\alpha)
        Y_\l (\alpha^\prime), \quad A_\l (\alpha) := \int_{S^2}
        A(\alpha^\prime, \alpha) \overline{Y_\l(\alpha^\prime)}
        d \alpha^\prime,
\eqno{(2.1)}
$$
 where $S^2$ is the unit sphere in $\R^3$,
$Y_\l (\alpha^\prime) = Y_{\l,m} (\alpha^\prime), -\l \leq m \leq 
\l$, are the normalized spherical harmonics,
summation over $m$ is understood in (2.1) and in (2.8) below.
Define the following algebraic variety:
$$
  M := \{ \theta : \theta \in \C^3, \theta \cdot \theta =1\}, \quad
        \theta \cdot w := \sum^3_{j=1} \theta_j w_j.
\eqno{(2.2)}
$$
This variety is non-compact, intersects $\R^3$ over $S^2$, and, given any
$\xi \in \R^3$, there exist (many) $\theta, \theta^\prime \in M$ such that
        $$\theta^\prime - \theta = \xi, \quad |\theta| \to \infty, \quad
        \theta, \theta^\prime \in M.  
\eqno{(2.3)}
$$
In particular, if one chooses the coordinate system in which
$\xi = te_3$, $t>0$, $e_3$ is the unit vector along the $x_3$-axis, then
the vectors
$$
        \theta^\prime = \frac{t}{2}e_3 + \zeta_2e_2 + \zeta_1e_1, \quad
        \theta = -\frac{t}{2}e_3 + \zeta_2e_2+\zeta_1 e_1, \quad
        \zeta^2_1 + \zeta^2_2 = 1-\frac{t^2}{4},
\eqno{(2.4)}
$$
satisfy (2.3) for any complex numbers $\zeta_1$ and $\zeta_2$ satisfying 
the last 
equation (2.4) and such that $|\zeta_1|^2+|\zeta_2|^2 \to \infty$. 
There are
infinitely many such $\zeta_1, \zeta_2 \in \C$.
Consider a subset $M^\prime \subset M$ consisting of the vectors
$\theta = (\sin \vartheta \cos \varphi, \sin \vartheta \sin \varphi, \cos 
\vartheta$
where $\vartheta$ and $\varphi$ run through the whole complex plane. 
Clearly
$\theta \in M$, but $M^\prime$ is a proper subset of $M$. Indeed, any
$\theta \in M$ with $\theta_3 \neq \pm 1$ is an element of $M^\prime$.
If $\theta_3 = \pm 1$, then $\cos \vartheta = \pm 1$, so
$\sin \vartheta = 0$ and one gets $\theta = (0,0, \pm 1) \in M^\prime$. 
However,
there are vectors $\theta = (\theta_1, \theta_2, 1) \in M$ which do not
belong to $M^\prime$. Such vectors one obtains choosing  
$\theta_1, \theta_2 \in \C$ such that $\theta^2_1 + \theta_2^2 = 0$.
There are infinitely many such vectors. The same is true for vectors
$(\theta_1, \theta_2, -1)$. Note that in (2.3) one can replace $M$ by
$M^\prime$ for any $ \xi \in \R^3$, $\xi \neq 2e_3$.

Let us state two estimates proved in \cite{R2}:
$$        \max_{\alpha \in S^2} \left| A_\l (\alpha) \right| \leq
        c \left(\frac{a}{\l}\right)^{\frac{1}{2}}
        \left(\frac{ae}{2\l}\right)^{\l +1},
\eqno{(2.5)}
$$        
and
$$        \left|Y_\l (\theta)\right| \leq \frac{1}{\sqrt{4 \pi}}
        \frac{e^{r |Im \theta|}}{|j_\l (r)|}, \quad
        \forall r > 0, \quad \theta \in M^\prime,
\eqno{(2.6)}
$$
where
$$        j_\l (r)  :=\left(\frac{\pi}{2r} \right)^{\frac{1}{2}}
        J_{\l + \frac{1}{2}} (r) = \frac{1}{2\sqrt{2}} \frac{1}{\l}
        \left(\frac{er}{2\l}\right)^\l [1 + o(1)] \hbox{\ as\ }
        \l \to \infty,
\eqno{(2.7)}
$$
and $J_\l(r)$ is the Bessel function regular at $r=0$.
Note that $Y_\l (\alpha^\prime)$,
defined above, admits a natural analytic continuation from $S^2$ to 
$M$
by taking $\vartheta$ and $\varphi$  to be arbitrary complex 
numbers.
The resulting $\theta^\prime \in M^\prime \subset M$.

The series (2.1) converges absolutely and uniformly on the sets
$S^2 \times  M_c$, where $M_c$ is any compact subset of $M$.

Fix any numbers $a_1$ and $b$, such that $a<a_1<b$. Let 
$||\cdot||$ denote the $L^2(a_1\leq |x|\leq b)$-norm. If $|x|\geq a,$
then the scattering solution is given analytically:
$$
 u(x, \alpha) = e^{i \alpha \cdot x} + \sum^\infty_{\l=0}
        A_\l (\alpha) Y_\l (\alpha^\prime) h_\l (r), \quad
        r:= |x| > a, \quad \alpha^\prime := \frac{x}{r},
\eqno{(2.8)}
$$
where $A_\l (\alpha)$ and $Y_\l(\alpha^\prime)$
are defined above,
$$h_\l(r) := e^{i \frac{\pi}{2} (\l+1)} \sqrt{\frac{\pi}{2r}}
H^{(1)}_{\l + \frac{1}{2}} (r),$$
$H^{(1)}_\l (r)$ is the Hankel function, and the normalizing factor is
chosen so that
$        h_\l(r) = \frac{e^{ir}}{r} [1 + o(1)] \hbox{\ as\ } r \to 
\infty.$
Define
$$
 \rho(x) := \rho(x;\nu) := e^{-i \theta \cdot x} \int_{S^2}
        u(x, \alpha) \nu(\alpha, \theta)d \alpha -1, \quad \nu \in 
L^2(S^2).
\eqno{(2.9)}
$$
Consider the minimization problem
$$        \Vert \rho \Vert = \inf := d(\theta),
\eqno{(2.10)} 
$$
where the infimum is taken over all $\nu \in L^2(S^2)$, and (2.3) holds.
        
It is proved in \cite{R2}  that
$$        d(\theta) \leq c |\theta|^{-1} \hbox{\ if\ } \theta \in M, \quad
        |\theta| \gg 1.
\eqno{(2.11)} 
$$
The symbol $|\theta| \gg 1$ means that $|\theta|$ is sufficiently large. 
The
constant $c>0$ in (2.11) depends on the norm $\Vert q \Vert_{L^2(B_a)}$ 
but
not on the potential $q(x)$ itself. 

An algorithm for computing a function
$\nu(\alpha, \theta)$, which can be used for inversion of the 
exact, fixed-energy, three-dimensional
scattering data, is as follows:

a) Find an approximate solution to (2.10) in the sense
$$        \Vert \rho (x, \nu) \Vert < 2 d(\theta),
\eqno{(2.12)}
$$
where in place of the factor 2 in (2.12) one could put any fixed constant 
greater
than 1.

b) Any such $\nu (\alpha, \theta)$ generates an estimate of
$\widetilde {q}(\xi)$ with the error $O\left(\frac{1}{|\theta|}\right)$,
$|\theta | \to \infty$.
This estimate is calculated by the formula
$$        \widehat{q} := -4\pi \int_{S^2} A(\theta^\prime, \alpha)  
        \nu (\alpha, \theta) d \alpha,
\eqno{(2.13)}
$$
where $\nu(\alpha, \theta) \in L^2(S^2)$ is any function satisfying
(2.12).

Our basic result is:
        
\begin{theorem}  Let (2.3) and (2.12) hold. Then
$$     \sup_{\xi \in \R^3}   \left|\widehat{q} - \widetilde{q}(\xi) 
\right| \leq
        \frac{c}{|\theta|}, \quad
        |\theta| \to \infty, \quad
\eqno{(2.14)}
$$
The constant $c>0$ in (2.14) depends 
on the norm of $q$, but not on a particular $q$.
\end{theorem}

In \cite{R2} and \cite{R7} an inversion algorithm is formulated also for
noisy data, and the error estimate for this algorithm is obtained.
 Let us describe these results.

Assume that the scattering data are given with some error: a function
$A_\delta (\alpha^\prime, \alpha)$ is given such that
$$
        \sup_{\alpha^\prime, \alpha \in S^2}
        \left|A(\alpha^\prime, \alpha) -
        A_\delta(\alpha^\prime, \alpha) \right| \leq \delta.
\eqno{(2.15)}
$$

We emphasize that $A_\delta (\alpha^\prime, \alpha)$ is not necessarily a
scattering amplitude corresponding to some potential, it is an arbitrary
 function in $L^\infty (S^2 \times S^2)$ satisfying (2.15). It is assumed 
that
the unknown function $A(\alpha^\prime, \alpha)$ is the scattering 
amplitude
corresponding to a $q \in Q$.
        
The problem is: {\it Find an algorithm for calculating 
$\widehat{q_\delta}$
such that}
$$        \sup_{\xi \in \R^3} \left|\widehat{q_\delta} - \widetilde{q} 
(\xi)
        \right| \leq \eta(\delta), \quad \eta(\delta) \to 0
        \hbox{\ as\ } \delta \to 0,
\eqno{(2.16)}
$$
{\it  and estimate
the rate at which $\eta(\delta)$ tends to zero.}

An algorithm for inversion of noisy data will now be described.

Let
$$
        N(\delta) := \left[ \frac{|\ln \delta|}{\ln|\ln \delta|} \right],  
\eqno{(2.17)}
$$
where $[x]$ is the integer nearest to $x>0$,
$$        \widehat{A}_\delta (\theta^\prime, \alpha) :=
        \sum^{N(\delta)}_{\l =0} A_{\delta \l} (\alpha) Y_\l
        (\theta^\prime), \quad
        A_{\delta \l} (\alpha) := \int_{S^2}
        A_\delta (\alpha^\prime, \alpha) \overline{Y_\l(\alpha^\prime)}
        d \alpha^\prime,
\eqno{(2.18)}
$$   
$$        u_\delta (x, \alpha) := e^{i \alpha \cdot x} +
        \sum^{N(\delta)}_{\l =0} A_{\delta \l}
        (\alpha) Y_\l (\alpha^\prime) h_\l (r),
\eqno{(2.19)}
$$   
$$        \rho_\delta (x; \nu) := e^{-i \theta \cdot x} \int_{S^2}
        u_\delta (x, \alpha) \nu (\alpha) d \alpha -1, \quad
        \theta \in M,
\eqno{(2.20)}   
$$
$$        \mu(\delta) := e^{-\gamma N(\delta)}, \quad
        \quad \gamma=\ln \frac {a_1}{a} >0,
\eqno{(2.21)}   
$$
$$
        a(\nu) := \Vert \nu \Vert_{L^2(S^2)}, \quad
        \kappa := |Im \theta|.
\eqno{(2.22)}   
$$
Consider the variational problem with constraints:
$$        |\theta| = \sup := \vartheta(\delta),
\eqno{(2.23)}   
$$
$$        |\theta| \left[ \Vert \rho_\delta (\nu) \Vert +
        a(\nu) e^{\kappa b} \mu(\delta) \right] \leq c,
        \quad \theta \in M, \quad
        |\theta| = \sup := \vartheta(\delta),
\eqno{(2.24)}
$$
the norm is defined above (2.8), and it is assumed that (2.3) holds,
where $\xi \in \R^3$ is an arbitrary fixed vector, $c>0$ is a sufficiently
large constant, and the supremum is taken over $\theta \in M$ and
$\nu \in L^2(S^2)$ under the constraint (2.24). By $c$ we denote various 
positive constants.

Given $\xi \in \R^3$ one can always find $\theta$ and $\theta^\prime$
such that (2.3) holds.
  We prove that $\vartheta(\delta) \to \infty$, more precisely:
$$        \vartheta(\delta) \geq c \frac{|\ln \delta|}{(\ln|\ln 
\delta|)^2},
        \quad \delta \to 0.
\eqno{(2.25)}
$$

Let the pair
$\theta(\delta)$ and $\nu_\delta (\alpha, \theta)$ be any approximate 
solution to
problem (2.23)-(2.24) in the sense that
$$        |\theta(\delta)| \geq \frac{\vartheta(\delta)}{2}.
\eqno{(2.26)}
$$
Calculate
$$        \widehat{q}_\delta := -4 \pi \int_{S^2} \widehat{A}_\delta
        (\theta^\prime, \alpha) \nu_\delta (\alpha, \theta) d \alpha.
\eqno{(2.27)}
$$

\begin{theorem}
      If (2.3) and (2.26) hold, then
$$        \sup_{\xi \in \R^3} \left| \widehat{q}_\delta -
        \widetilde{q}(\xi) \right| \leq c
        \frac{(\ln |\ln \delta|)^2}{|\ln \delta|} \hbox{\ as\ }
        \delta \to 0,
\eqno{(2.28)}
$$
where $c>0$ is a constant depending on the norm of $q$.
\end{theorem}   
 In [8] estimates (2.14) and (2.28) were formulated with the supremum
taken over an arbitrary large but fixed ball of radius $\xi_0$. Here these
estimates are improved: $\xi_0=\infty$. The key point is: the constant
$c>0$ in the estimate (2.11) does not depend on $\xi$.

{\bf Remark.} In \cite{R8} (see also \cite{R9} and \cite{R7}) an analysis 
of the approach to ISP, based on the recovery of the DN 
(Dirichle-to-Neumann) map from the fixed-energy scattering data, is given.
This approach is discussed in Section 3.

The basic numerical difficulty of the approach described in Theorems 
2.1 
and 2.2 comes from solving problems (2.10) for exact data, and 
problem (2.23)-(2.24)  for noisy data. Solving (2.10) amounts to
finding a global minimizer of a quadratic form of the variables $c_\ell$, 
if one 
takes $\nu$ in (2.9) as a linear combination of the spherical harmonics:
$\nu=\sum_{\ell =0}^L c_\ell  Y_\ell (\alpha)$. If one uses
the necessary condition for a minimizer of a quadratic form, that is, a 
linear system, then the matrix of this system is ill-conditioned
for large $L$. This causes the main difficulty in the numerical 
solution of (2.10). On the other hand, there are methods for global 
minimization of the quadratic functionals, based on the gradient descent,
which may be more efficient than using the above necessary condition.
        
\section{Discussion of the inversion method which uses the DN map}

In [14] the following inversion method is discussed:

$$
\tilde q (\xi)=\lim_{|\theta|\to \infty} \int_S\exp
(-i\theta'\cdot s)(\Lambda-\Lambda_0)\psi ds , \eqno {(3.1)}
$$
where (2.3) is assumed, $\Lambda$ is the Dirichlet-to-Neumann (DN) map,
$\psi$ is 
found from the equation:
$$
\psi (s)=\psi_0(s)-\int_S G(s-t)B\psi dt,\quad B:=\Lambda-\Lambda_0 .\eqno
{\eqno (3.2)}
$$
and $G$ is defined by the formula:
$$ 
G(x)=\exp (i\theta\cdot x)\frac 1{(2\pi)^3}\int_{\mathbb R^3}
\frac{\exp (i\xi\cdot x)d\xi}{\xi^2+2\xi\cdot\theta}.\eqno {(3.3)}
$$
The DN map is constructed from the 
fixed-energy scattering data $A(\alpha', \alpha)$
by the method of [14] (see also [8]).

Namely, given $A(\alpha', \alpha)$ for all $\alpha', \alpha \in S^2$,
one finds $\Lambda$ using the following steps.

Let $f\in H^{3/2}(S)$ be given, $S$ is a sphere of radius $a$ centered 
at the origin, $f_\l$ are its Fourier coefficients
in the basis of the spherical harmonics, 
$$
w=\sum^\infty_{l=0}f_lY_l(x^0)\frac{h_l(r)}{h_l(a)},\quad r\geq a,
\quad x^0:=\frac xr, \quad r:=|x|.\eqno {(3.4)}
$$
Let 
$$
w=\int_S g(x,s)\sigma (s)ds ,\eqno {(3.5)}
$$  
where $\sigma$ is some function, which we find below, and $g$ is the 
Green function (resolvent kernel) of the Schroedinger operator, 
satisfying the radiation condition at infinity. 
Then
$$
w^+_N=w^-_N+\sigma ,\eqno {(3.6)}
$$
where $N$ is the outer normal to $S$, so $N$ is directed along the 
radius-vector. 
We require $w=f$ on $S$. Then $w$ is given by (3.4) in the exterior of 
$S$, and 
$$
w_N^-=\sum^\infty_{l=0}f_lY_l(x^0)\frac{h^\prime_l(a)}{h_l(a)}.\eqno {(3.7)}
$$
By formulas (3.6) and (3.7), finding $\Lambda$ is equivalent to finding 
$\sigma$. By (3.5), asymptotics of $w$ as $r:=|x|\to \infty$, 
$x/|x|:=x^0,$ is (cf [8, p.67]):
$$
w=\frac {e^{ir}}{r} \frac {u(y,-x^0)}{4\pi} +o(\frac 1 r),
\eqno {(3.8)}
$$
where $u$ is the scattering solution,
$$
u(y,-x^0)=e^{-ix^0 \cdot y}+\sum_{\l=0}^\infty A_\l 
(-x^0)Y_\l(y^0)h_\l(|y|). \eqno {(3.9)}
$$
From (3.4), (3.8) and (3.9) one gets an equation for finding $\sigma$
([14, eq. (23)], see also [8, p. 199]):
$$
\frac{f_l}{h_l(a)}=\frac 1{4\pi}\int_S ds\sigma (s)\left (u(s,-\beta),
Y_l(\beta)\right )_{L^2(S^2)} ,\eqno {(3.10)}
$$
which can be written as a linear system:
$$
\frac{4\pi f_l}{h_l(a)}=a^2(-1)^l\sum^\infty_{l^\prime
=0}\sigma_{l^\prime}
[4\pi i^lj_l(a)\delta_{ll^\prime}+A_{l^\prime l}h_{l^\prime}(a)] ,\eqno
{\eqno (3.11)}
$$
for the Fourier coefficients $\sigma_\l$ of $\sigma$. The coefficients
$$A_{l^\prime l}:=((A(\alpha', \alpha), Y_\l(\alpha'))_{L^2(S^2)}, Y_\l 
(\alpha))_{L^2(S^2)} 
$$ 
are the Fourier coefficients of the scattering
amplitude. Problems (3.10) and (3.11) are very ill-posed 
(see [14] for details).

This approach faces heavy difficulties:

1) The construction of the DN map from the scattering data is a very
ill-posed problem,

2) The construction of the potential
from the DN map is a very difficult problem numerically,  because
one has to solve a Fredholm-type integral equation
( equation (3.2) ) whose kernel contains $G$, defined in (3.3). This $G$ 
is a
tempered distribution, and it is
very difficult to compute it,

3) One has to calculate a limit
of an integral whose integrand grows exponentially to infinity
if a factor in the integrand is not known exactly. The solution of
equation (3.2)  is one of the factors in the integrand.
It cannot be known exactly in practice because it cannot be
calculated with arbitrary accuracy even if the scattering data are known
exactly. Therefore the limit in formula (3.1)  cannot be
calculated accurately.

No error estimates are obtained for this approach.
 
In contrast, in Ramm's method, there is no need to 
compute $G$, to solve equation (3.2), to calculate the DN map from
the scattering data, and to compute the limit (3.1).
The basic difficulty in Ramm's inversion method for exact data is to 
minimize a quadratic form (2.10), and for noisy data to solve optimization 
problem (2.23)-(2.24). The error estimates are obtained for the Ramm's
method.

\end{document}